\documentclass[a4paper]{article}
\usepackage[T1]{fontenc}
\usepackage[latin1]{inputenc}
\usepackage{amsmath,amsthm,amsfonts,amssymb, euscript, stmaryrd, graphicx,superbox}
\usepackage{algorithm}
\usepackage{algorithmic}
\usepackage{color}
\usepackage{listings}
\usepackage{pstricks}
\usepackage{psfrag}
\usepackage{amsthm}
\usepackage{fullpage}
\definecolor{colKeys}{rgb}{0,0,1}
\definecolor{unidentified}{rgb}{0,0,0}
\definecolor{concealments}{rgb}{0,0.5,0}
\definecolor{colString}{rgb}{0.6,0.1,0.1}

\definecolor{Grey1}{rgb}{0.5,0.5,0.5}
\definecolor{Grey2}{rgb}{0.8,0.8,0.8}
\definecolor{Grey3}{rgb}{0.9,0.9,0.9}
\definecolor{Grey4}{rgb}{0.95,0.95,0.95}
\definecolor{ACORRIGER}{rgb}{0.9,0,0}

\newcommand{\footnoteremember}[2]{
 \footnote{#2}
 \newcounter{#1}
 \setcounter{#1}{\value{footnote}}
}

\author
{
Sylvain \textsc{Corlay}\footnote{Natixis, Equity Derivatives and Arbitrage. E-mail: sylvain.corlay@natixis.com. The authors would like to thank the members of Natixis equity derivatives quantitative R\&D team for fruitful discussions.}
\footnoteremember{myfootnote}{Laboratoire de Probabilités et Modèles Aléatoires, UMR 7599, Université Paris 6, case 188, 4, pl. Jussieu, F-75252 Paris Cedex 5, France. }\\
}

\title{The Nyström method for functional quantization with an application to the fractional Brownian motion}
\date{September 6, 2010}

\newtheorem{theo}{Theorem}[section]
\newtheorem{prop}[theo]{Proposition}

\newcommand{\1}{\textbf{1}}
\newcommand{\N}{\mathbb{N}}

\newcommand{\R}{\mathbb{R}}

\newcommand{\E}{\mathbb{E}}

\newcommand{\PP}{\mathbb{P}}
\newcommand{\diag}{\operatorname{diag}}

\newcommand{\card}{\operatorname{card}}

\newcommand{\sspan}{\operatorname{span}}
\newcommand{\cov}{\operatorname{cov}}
\newcommand{\var}{\operatorname{Var}}

\newcommand{\ui}{{\underline{i}}}

\def\keywordname{{\bf Keywords:}} 
\newcommand{\keywords}[1]{\par\addvspace\baselineskip\noindent\keywordname\enspace\ignorespaces#1}

\graphicspath{{./}{figures/}}
\begin{document}
\lstset{language=c++}

%
%		TITLE
%

\maketitle

%
%		ABSTRACT
%

\begin{abstract}
\par In this article, the so-called "Nyström method" is tested to compute optimal quantizers of Gaussian processes. In particular, we derive the optimal quantization of the fractional Brownian motion by approximating the first terms of its Karhunen-Loève decomposition. 
\par A numerical test of the "functional stratification" variance reduction algorithm is performed with the fractional Brownian motion.
\end{abstract}

\pagebreak
\section*{Introduction}\label{sec:introduction}
\par Let $(\Omega,\mathcal{A},\PP)$ be probability space, and $E$ a reflexive Banach space. The norm on $E$ is denoted $|\cdot |$. 
\par The quantization of a random variable $X$, taking its values in $E$ consists in its approximation by a random variable $Y$ taking finitely many values. The resulting error of this discretization is the $L^p$ norm of $|X-Y|$. Minimizing this error, with a fixed maximum cardinal of $Y(\Omega)$ yields the following minimization problem:
\begin{equation}
\min \left\{ \| X-Y \|_p, Y:\Omega \to E \textrm{ measurable }, \card(Y(\Omega)) \leq N \right\}.
\end{equation}
\par \noindent This problem, was first considered for signal transmission and compression issues. More recently, quantization has been introduced in numerical probability, to devise quadrature methods \cite{PagesIntegVectorQuant}, solving multi-dimensional stochastic control problems \cite{BallyPagesPrintemsAmerican1} and for variance reduction \cite{CorlayPagesStratification}. Since the $2000$'s, the infinite dimensional setting has been investigated from both theoretical an numerical viewpoint, especially in the quadratic case \cite{LuschgyPagesFunctional3}. One elementary property of a $L^2$ optimal quantizer is the stationarity: $\E[X|Y] = Y.$
\par If $X$ is a bi-measurable stochastic process on $[0,T]$ verifying $\int_0^T \E[|X_t|^2] dt < \infty$, it can be considered as a random variable valued in the Hilbert space $H = L^2([0,T])$. In \cite{LuschgyPagesFunctional3}, it is shown that in the centered Gaussian case, linear subspaces $U$ of $H$ spanned by $N$-stationary quantizers correspond to principal components of $X$, in other words, are spanned by eigenvectors of the covariance operator of $X$. Thus, the quantization consists first in exploiting its Karhunen-Loève decomposition $\left(e_n^X, \lambda_n^X\right)_{n \geq 1}$. 
\par If $d^X(N)$ is the dimension of the subspace of $L^2([0,T])$ spanned by $Y(\Omega)$, the quantization error $e_N(Y)$ writes
\begin{equation}\label{eq:distorsion_representation}
e_N(X)^2 = \sum\limits_{j \geq m+1} \lambda_j^X + e_N\left( \bigotimes\limits_{j=1}^m \mathcal{N}(0,\lambda_j^X)\right)^2 \textrm{ for } m \geq d_N(X).
\end{equation}
\begin{equation}
e_N(X)^2 < \sum\limits_{j \geq m+1} \lambda_j^X + e_N\left( \bigotimes\limits_{j=1}^m \mathcal{N}(0,\lambda_j^X)\right)^2 \textrm{ for } 1 \leq m < d_N(X).
\end{equation}
\par The decomposition is first truncated at a fixed order $m$ and then the $\mathbb{R}^m$-value Gaussian vector constituted of the $m$ first coordinates of the process on its Karhunen-Loève decomposition is quantized. To reach optimal quantization, we have both to determine the optimal rank of truncation $d^X(N)$ (the quantization dimension) and to determine the optimal $d^X(N)$-dimensional Gaussian quantizer corresponding to the first coordinates, $\bigotimes\limits_{j=1}^{d^X(N)} \mathcal{N}(0,\lambda^X_j)$. Usual examples of such processes are the standard Brownian motion on $[0,T]$, the standard Brownian bridge on $[0,T]$, the fractional Brownian motion and the fractional Ornstein-Uhlenbeck process.
\par We can also choose to use a product quantization of $\bigotimes\limits_{j=1}^m \mathcal{N}(0,\lambda^X_j)$. The product quantization is the cartesian product of the optimal quantizers of the standard one-dimensional Gaussian distributions $\mathcal{N}\left(0,\lambda^X_i\right)_{1 \leq i \leq d^X(N)}$. In the case of independent marginals, this yields a stationary quantizer. One advantage of this method is that the one-dimensional Gaussian quantization is a fast procedure. Newton-Raphson methods converge very fast to the optimal quantization (see \cite{PagesGaussianQuantization}). Moreover, a sharply optimized database of quantizers of standard univariate and multivariate Gaussian distributions is available on the web site {\verb www.quantize.maths-fi.com } \cite{WebSiteGaussian} for download. Still, we have to determine quantization size on each dimension to obtain optimal product quantization. In this case, the minimization of the distorsion (\ref{eq:distorsion_representation}) comes to:
\begin{equation}\label{eq:product_error_minimization}
\min\left\{\sum\limits_{n=1}^d \lambda_n^X \min\limits_{\R^{N_n}}\|\xi-\tilde{\xi}^{(N_n)} \|_2^2 + \sum\limits_{n \geq d+1} \lambda_n^X, N_1 \times \cdots \times N_d \leq N, d \geq 1 \right\}.
\end{equation}
A solution of (\ref{eq:product_error_minimization}) is called an optimal K-L product quantizer. This problem can be solved by the "blind optimization procedure", which consists in computing the criterium for every possible decomposition $N_1 \times \cdots \times N_d$ with $N_1 \geq \cdots \geq N_d$. The result of this procedure can be kept off-line for a future use. Optimal decompositions for a wide range of values of $N$ for both Brownian motion and Brownian bridge are available on the web site {\verb www.quantize.maths-fi.com } \cite{WebSiteGaussian}. 
\par In \cite{LuschgyPagesFunctional3}, the rate of convergence to zero of the quantization error is investigated. A complete solution is provided for the case of Gaussian processes with regular varying eigenvalues. Rates of convergence are available for the above cited examples of Gaussian processes. The asymptotic of the quantization dimension $d^X(N)$ are investigated in \cite{LuschgyPagesFunctional2}. The following theorem combines these results:
\begin{theo}[Functional quantization asymptotics]\label{thm:quantization_asymptotics}
Let $X$ be a centered Gaussian process on $[0,T]$ with Karhunen-Loève system $(e_n^X,\lambda_n^X)_{n \geq 1}$. Let $(Y_N)_{N\geq 1}$ be a sequence of quadratic optimal $N-$quantizers for $X$. We assume that
$$
\lambda_n^X \sim \frac{\kappa}{n^b} \textrm{ as } n \to \infty \hspace{5mm} (b>1).
$$
We have:
\begin{itemize}
\item $\sspan(Y_N(\Omega)) = \sspan\left\{e_1^X, \cdots, e_{d^X(N)}^X \right\}$ and $d^X(N) = \Omega(\log N )$.
\item $e_N(X) =\|X-Y_N \|_2 \sim \sqrt{\kappa}\sqrt{b^b (b-1)^{-1}} (2 \log N)^{-\frac{b-1}{2}}$
\end{itemize}
\end{theo}
\par \noindent A conjecture is $d^X(N) \sim \frac{2}{b} \log(N)$.
\vspace{5mm}
\par \noindent It is shown in \cite{LuschgyPagesFunctional3} that the Karhunen-Loève eigenvalues of the fractional Brownian motion, $(\lambda^{B^H}_n)_{n \geq 1}$ verify 
$$
\lambda^{B^H}_n \sim \frac{1}{n^{2H+1}} \textrm{ as } n \to \infty,
$$
thus the fractional Brownian motion satifies the hypothesis of theorem \ref{thm:quantization_asymptotics}.
\par \noindent In a constructive viewpoint, the numerical computation of the optimal quantization or the optimal product quantization requires a numerical evaluation of the Karhunen-Loève eigenfunctions and eigenvalues, at least the very first terms. (As seen in theorem \ref{thm:quantization_asymptotics}, the quantization dimension of usual Gaussian processes increases asymptotically as the logarithm of the size of the quantizer, so it is most likely that it is small. For instance, the quantization dimension $d^W(N)$ of the Brownian motion with $N=10000$ is $9$.) The Karhunen-Loève decomposition of some usual Gaussian processes have a closed-form expression. It is the case of the standard Brownian motion, the Brownian bridge and the Ornstein-Uhlenbeck process. (The special case of the Ornstein-Uhlenbeck process is derived in \cite{CorlayPagesStratification}).
{
\small
\begin{enumerate}
	\item The Brownian motion $(W_t)_{t \in [0,T]}$,
	\begin{equation}\label{eq:brownian_motion_kl}
	e_n^W(t) := \sqrt{\frac{2}{T}} \sin \left( \pi (n-1/2) \frac{t}{T}\right), \hspace{8mm} \lambda_n^W:= \left(\frac{T}{\pi (n-1/2)} \right)^2, \hspace{4mm} n \geq 1.
	\end{equation}
	\item The Brownian bridge on $[0,T]$, 
	\begin{equation}\label{eq:brownian_bridge_kl}
	e_n^B(t) := \sqrt{\frac{2}{T}} \sin \left( \pi n \frac{t}{T}\right), \hspace{8mm} \lambda_n^B:= \left(\frac{T}{\pi n} \right)^2, \hspace{4mm} n \geq 1.
	\end{equation}
	\item The Ornstein-Uhlenbeck process on $[0,T]$, starting from $0$, defined by the SDE $dr_t = \theta (mu - r_t)dt + \sigma dW_t$, with $\sigma \geq 0$, $\theta>0 $ and $W$ a standard Brownian motion on $[0,T]$.
	\begin{equation}
	e_n^{OU}(t) := \left( \frac{1}{\sqrt{\frac{T}{2}-\frac{\sin(2\omega_{\lambda_n}T)}{4\omega_{\lambda_n}}}} \right) \sin(\omega_{\lambda_n} t), \hspace{8mm} \lambda_n^{OU}:= \frac{\sigma^2}{\omega_{\lambda_n}^2 + \theta^2}, \hspace{4mm} n \geq 1,
	\end{equation}
	where $\omega_{\lambda_n}$ are the (sorted) strictly positive solutions of the equation 
	$$
	\theta \sin(\omega_{\lambda_n} T) + \omega_{\lambda_n} \cos(\omega_{\lambda_n} T) = 0.
	$$
	\item The stationary Ornstein-Uhlenbeck process on $[0,T]$, defined by the same SDE with $r_0 \sim \mathcal{N}(0,\sigma_0)$.
	\begin{equation}
	e_n^{OU}(t) := C_n \left(\omega_{\lambda_n} \cos(\omega_{\lambda_n} t) + \theta \sin( \omega_{\lambda_n}t) \right), \hspace{8mm} \lambda_n^{OU}:= \frac{\sigma^2}{\omega_{\lambda_n}^2 + \theta^2}, \hspace{4mm} n \geq 1,
	\end{equation}
	where $\omega_{\lambda_n}$ are the (sorted) strictly positive solutions of the equation 
	$$
	2 \theta \omega \cos(\omega_{\lambda_n} T) + (\theta^2-\omega_{\lambda_n}^2) \sin(\omega_{\lambda_n} T) = 0,
	$$
	and
	$$
	\frac{1}{C_n^2} = \frac{\theta}{2} \left(1-\cos(2 \omega_{\lambda_n}T )\right) + \frac{\omega_{\lambda_n}}{2} \left( T + \frac{\sin(2 \omega_{\lambda_n}T)}{2 \omega_{\lambda_n}} \right) +\frac{\theta^2}{2} \left( T - \frac{\sin(2 \omega_{\lambda_n} T)}{2 \omega_{\lambda_n}} \right).
	$$
\end{enumerate}
}
\par In a more general setting, we do not have a closed-form expression for the Karhunen-Loève decomposition. For instance, as far as we know, the K-L expansion of the fractional Brownian motion is not known. Hence, a numerical method to evaluate first Karhunen-Loève eigenfunctions is the "missing link" on the path to the constructive optimal quantization of more Gaussian processes. 
\par However, we can derive rate-optimal quantization of Gaussian processes using other series expansions as proposed by Luschgy and Pages in \cite{ProductQuantLuschyPages,LuschyPagesParseval}. In this setting, the case of the fractional Brownian motion can be derived using a rate-optimal series expansion proved by Dzhaparidze and van Zanten in \cite{DzhaparidzevanZanten1,DzhaparidzevanZanten2}. Other constructive approaches for functional quantization are proposed by Wilbertz in \cite{WilbertzPHD}.
\par In this article, we experiment the so-called "Nyström method" \cite{AtkinsonMonograph,DelvesMohamedIntegral,NumericalRecipes} for approximating the solution of the functional eigenvalue problem which defines the Karhunen-Loève decomposition. First, we compare the result of the the numerical method with the closed-forms available for the Brownian motion, the Brownian bridge and the Ornstein-Uhlenbeck process. Then, the special case of the functional quantization of the fractional Brownian motion is handled. 
\par Functional quantization of Gaussian processes have numerous applications in numerical probability. In \cite{CorlayPagesStratification}, a variance reduction method based on the functional quantization of a Gaussian process was proposed. This method can be seen as a "Guided Monte-Carlo simulation" (see figure \ref{fig:fractional_brownian_motion_conditional_simulation}). Still, it was only applicable with Gaussian processes for which we could have a numerical evaluation of the Karhunen-Loève eigenfunctions. Such a variance reduction method would be of high interest in Monte-Carlo simulations implying the fractional Brownian motion because its simulation schemes have a high complexity. 
\par Subsequently, we test this "functional stratification" variance reduction algorithm in option pricing problems within the fractional Brownian motion's counterpart of the classical Black and Scholes model. First, the case of a Vanilla option is benchmarked with the closed-form expression available in this case. Then the case of discrete barrier options is tested.

\section{The Nyström method}\label{seq:nystrom_introduction}
\par Let $X$ be a bi-measurable Gaussian stochastic process on $[0,T]$ defined on the probability space $(\Omega,\mathcal{A},\PP)$. We assume that $\int_{[0,T]} \E[X_s^2]ds<\infty$. Let us denote 
$\Gamma_X(t,s)$ the covariance function of $X$ defined by $\Gamma_X(t,s) = \cov(X_t,X_s)$. The covariance operator $C_X$ of $X$ is defined by $C_X f = \int_{[0,T]} \Gamma_X(\cdot,s) f(s) dt.$ It is a symmetric positive trace class operator on $L^2[0,T]$. The Karhunen-Loève basis associated with $X$, denoted $(e^X_n)_{n \geq 1}$ is the Hilbert basis of $L^2[0,T]$ constituted with eigenvectors of $C_X$ with decreasing eigenvalues. Now, we aim to solve numerically the eigenvalue problem
\begin{equation}\label{eq:KL_eigen_equation}
\int_0^T \Gamma_X(\cdot,s) f_k(s) d s = \lambda_k f_k, \hspace{5mm} k \geq 1.
\end{equation}
The Nyström method requires the choice of some quadrature rule $\int_0^T f(s) ds \sim \sum\limits_{i=1}^n w_j f(s_j).$ $(w_j)_{1 \leq j \leq n}$ is the sequence of the weights of the quadrature rule, while $(s_j)_{1 \leq j \leq n}$ are the abscissas where $f$ is evaluated. If we introduce this quadrature rule in equation (\ref{eq:KL_eigen_equation}), we get
\begin{equation}\label{eq:discretized_KL_eigen}
\sum\limits_{j=1}^n w_j \Gamma_X(t,s_j) f_k(s_j) = \lambda_k f_k(t) \hspace{8mm} t \in [0,T].
\end{equation}
Evaluating equation (\ref{eq:discretized_KL_eigen}) at the quadrature points yields 
\begin{equation}
\sum\limits_{j=1}^n w_j \Gamma_X(t_i,s_j) f_ k(s_j) = \lambda_k f_k(t_i) \hspace{8mm} i \in \{1,\cdots,n\}.
\end{equation}
Let $f$ be the vector {\small$\left(\begin{array}{cc}f_k(t_1)\\ \vdots \\ f_k(t_{n}) \end{array}\right)$}, $((K_{ij}))_{1 \leq i,j \leq n}$ the matrix $\left(\left(\Gamma_X(t_i,s_j)\right)\right)_{1 \leq i,j \leq n}$, $\lambda = (\textrm{diag}(\lambda_k))_{k=1\cdots n}$ and define $\tilde{K}_{ij} = K_{ij}w_j$. Then the eigenvalue problem becomes 
\begin{equation}\label{eq:matrix_eigen}
\tilde{K} f = \lambda f.
\end{equation}
\par Hence, within this approximation, the functional eigenvalue problem turns into a matrix eigenvalue problem. As $K$ is a covariance matrix, it is symmetric. However, since the weights are not equal for most quadrature rules, the matrix $\tilde{K}$ is not symmetric. As outlined in \cite{NumericalRecipes}, numerical methods for matrix orthogonalization are much simpler in the symmetric case. As a consequence, we should restore the symmetry if possible. The method proposed in \cite{NumericalRecipes} is the following:
\par We define the diagonal matrix $D = \textrm{diag}(w_j)$ and its square root $D^{1/2} = \textrm{diag}(\sqrt{w_j})$. Then equation (\ref{eq:matrix_eigen}) becomes
\begin{equation}\label{eq:reweighted_equation}
K \cdot D \cdot f = \lambda f.
\end{equation}
Multiplying by $D^{1/2}$, we get
\begin{equation}\label{eq:reweighted_equation_2}
\left( D^{1/2} \cdot K \cdot D^{1/2} \right) \cdot h = \lambda h, \hspace{3mm} \textrm{where} \hspace{3mm} h = D^{1/2} \cdot f.
\end{equation}
\par Equation (\ref{eq:reweighted_equation_2}) is now in the form of a symmetric eigenvalue problem. For square-integrable kernels (we stand in this case), this provides a good approximation of the $n$ highest eigenvalues. 
\subsection{Choice of the quadrature method}
\par Classical numerical methods for real symmetric matrix diagonalization are 
\begin{itemize}
\item The Jacobi transformation for symmetric diagonalization.
\item A tridiagonalization (by Givens or Householder reduction) followed by a QL algorithm with implicit shifts.
\end{itemize}
All these numerical methods have a $O(n^3)$ complexity. As a consequence, the natural choice for the quadrature method would be the highest order possible (A high order Bode's formula, or a Gaussian quadrature). 
\par However as pointed out in \cite{GuoqiangExtrapolation}, the Nyström method associated with the trapezoidal integration rule admits an asymptotic error expansion in even powers of the step sizes as soon as the covariance function is differentiable (or continuous and piecewise differentiable).  As a consequence, instead of using the high order integration rule, we prefer to use a Richardson-Romberg extrapolation on the result of the whole procedure with the trapezoidal quadrature formula. We could reach an accuracy which approaches the machine roundoff error on the first eigenvalues when we benchmark this method on the Brownian motion, the Brownian bridge or the Ornstein-Uhlenbeck process. Another argument for the trapezoidal rule is that we encountered some small instabilities on the eigenfunction evaluation when using higher order schemes. 
\subsection{Choice of the interpolation method}
\par \noindent The natural choice is to use equation (\ref{eq:discretized_KL_eigen}) as an interpolation method for evaluating $f_k$,
\begin{equation}\label{eq:good_interpolation}
f_k(t) = \frac{1}{\lambda_k} \sum\limits_{j=1}^n w_j K(t,s_j) f_k(s_j).
\end{equation}
\par \noindent The same Richardson-Romberg extrapolation can be performed between the values of $\sum\limits_{j=1}^n w_j K(t,s_j) f_k(s_j)$ with the different orders $n$ to compute this integral. The result is then divided by the extrapolated value of $\lambda_k$.
\par \noindent \textbf{A remark on the interpolation method}
\par \noindent One purpose of the quantization of a Gaussian process $X$, is to perform a quantization of a diffusion with respect $X$, as soon as such a stochastic integral can be defined. We can obtain a quantizer of the diffusion by inserting the quantizer of the Gaussian process in the diffusion equation written in the Stratonovich sense. The most accomplished study on this subject is \cite{PagesSellamiSDE}. In this case, we may also need a numerical approximation of the time-derivative of the eigenfunction in the Karhunen-Loève decomposition. This work is mostly specific to the Brownian motion but main results remain valid for continuous semi-martingales that satisfy the Kolmogorov criteria as the Brownian bridge and Ornstein-Uhlenbeck processes. 
\par Still, a future work could be to extend these results to diffusions with respect to the fractional Brownian motion and other related processes. If $\Gamma_X$ is (weakly) differentiable, a natural evaluation method for the derivative would be $f'_k(t) = \frac{1}{\lambda_k} \sum\limits_{j=1}^n w_j \partial_1 \Gamma_X(t,s_j) f_k(s_j).$
\par One problem is that this method yields an irregular derivative. For example, this yields a piecewise constant derivative in the case of the Brownian motion. This causes instabilities problems when using Runge-Kutta integration methods for ordinary differential equations, which rely on the regularity of the considered Cauchy problem. 
\par As a consequence, a more regular interpolation method can give more satisfactory results when dealing with diffusions. (Spline or rational interpolation methods for instance.)
\section{Benchmark on known Karhunen-Loève expansions}
\par In this section, we compare the numerical results obtained with the Nyström methods in cases where we have closed-form expression of the Karhunen-Loève expansion. The multi-steps Richardson-Romberg extrapolation consists in using the asymptotic error estimate of the method
$$
V = u_n + \frac{K_1}{n^2} + \frac{K_2}{n^4} + \cdots + O\left(\frac{1}{n^{2p}}\right).
$$
Writing this expression for $p$ different values of $n$ allows us to solve a $p \times p$ linear system to nullify the $p-1$ first orders of convergence. The three-steps Richardson-Romberg extrapolation with $n=p$, $n=l$ and $n=k$ gives the following solution : 
$$
\frac{U_k k^4 (m^2-l^2) + U_l l^4 (k^2 - m^2) + U_m m^4 (l^2 - k^2) }{(m^2-l^2)(l^2 m^2+k^4-m^2 k^2-l^2 k^2)}.
$$
This result is naturally invariant by any permutation of the coefficients $(k, m, l)$. We experienced less accurate results when using higher order Richardson-Romberg extrapolation, so we will settle for a three-steps extrapolation.
\subsection{Eigenvalues accuracy}
\par In tables \ref{tab:brownnian_motion_eigenvalues_equi} and \ref{tab:brownnian_bridge_eigenvalues_equi}, Karhunen-Loève eigenvalues of the Brownian motion and of the Brownian bridge on $[0,1]$ are reported. Table \ref{tab:stationary_ornstein_eigenvalues_equi} deals with the stationary Ornstein-Uhlenbeck on $[0,1]$ defined by the SDE
\begin{equation}\label{eq:sde_ornstein_stationary}
dr_t = -\theta r_t dt + \sigma dW_t, \hspace{10mm} r_0 \sim \mathcal{N}\left(0,\frac{1}{2}\right).
\end{equation}
First column gives the theoretical value given by the closed-form. Following columns give the value computed with the Nyström method with a regular step size with $25$, $50$ and $100$ points. Last column gives the absolute error of a $3$ steps Richardson-Romberg extrapolation method between $n=25$, $n=50$ and $n=100$. 
\begin{figure}[!ht]
\begin{center}
{
\hspace*{-3mm}
\begin{tabular}{|c|c|c|c|c|}
\hline
				& Trapezoidal	& Trapezoidal	& Trapezoidal	& Trapezoidal Nyström\\
Closed-form		&	Nyström		&	Nyström 	&	Nyström 	& $25-50-100$ Richardson-Romberg\\
				& $25$ points	& $50$ points	& $100$ points 	& absolute error\\
\hline \hline
$0.405284735$	& $0.405418094$		& $0.405318070$		& $0.405293068$		& $6.3727$e$-14$\\
\hline
$0.0450316372$	& $0.0451652077$	& $0.0450649853$	& $0.0450399714$	& $5.2269$e$-12$\\
\hline
$0.0162113894$	& $0.0163453833$	& $0.0162447639$	& $0.0162197259$	& $4.0448$e$-11$\\
\hline
$0.00827111703$	& $0.00840574996$	& $0.00830453112$	& $0.00827945541$	& $1.5607$e$-10$\\
\hline
$0.00500351524$	& $0.00513900777$	& $0.00503698224$	& $0.00501185691$	& $4.2896$e$-10$\\
\hline
\end{tabular}
\caption{Record of the $5$ highest eigenvalues of the Karhunen-Loève decomposition of the Brownian motion. }
\label{tab:brownnian_motion_eigenvalues_equi}
}
\end{center}
\end{figure}

\begin{figure}[!ht]
\begin{center}
{
\hspace*{-3mm}
\begin{tabular}{|c|c|c|c|c|}
\hline
				& Trapezoidal	& Trapezoidal	& Trapezoidal	& Trapezoidal Nyström\\
Closed-form 	&	Nyström		&	Nyström 	&	Nyström 	& $25-50-100$ Richardson-Romberg\\
				& $25$ points 	& $50$ points	& $100$ points	& absolute error\\
\hline \hline
$0.101321184$	& $0.101454622$		& $0.101354524$		& $0.101329517$		& $1.0314$e$-12$\\
\hline
$0.0253302959$	& $0.0254640514$	& $0.0253636556$	& $0.0253386309$	& $1.6540$e$-11$\\
\hline
$0.0112579093$	& $0.0113921955$	& $0.0112913019$	& $0.0112662463$	& $8.4041$e$-11$\\
\hline
$0.00633257398$	& $0.00646760876$	& $0.00636601285$	& $0.00634091389$	& $2.6697$e$-10$\\
\hline
$0.00405284735$	& $0.00418885438$	& $0.00408634582$	& $0.00406119097$	& $6.5608$e$-10$\\
\hline
\end{tabular}
\caption{Record of the $5$ highest eigenvalues of the Karhunen-Loève decomposition of the Brownian bridge. }
\label{tab:brownnian_bridge_eigenvalues_equi}
}
\end{center}
\end{figure}

\begin{figure}[!ht]
\begin{center}
{
\hspace*{-3mm}
\begin{tabular}{|c|c|c|c|c|}
\hline
				& Trapezoidal	& Trapezoidal	& Trapezoidal	& Trapezoidal Nyström\\
Closed-form 	&	Nyström		&	Nyström 	&	Nyström 	& $25-50-100$ Richardson-Romberg\\
				& $25$ points 	& $50$ points	& $100$ points	& absolute error\\
\hline \hline
$0.369405405$	& $0.369395812$		& $0.369403011$		& $0.369404807$	& $2.7645$e$-13$\\
\hline
$0.0690018877$	& $0.0690750142$	& $0.0690201680$	& $0.0690064577$	& $2.0265$e$-12$\\
\hline
$0.0225442436$	& $0.0226553722$	& $0.0225719721$	& $0.0022551172$	& $5.3713$e$-12$\\
\hline
$0.0106644656$	& $0.0107875835$	& $0.0106950942$	& $0.0106721134$	& $5.8762$e$-11$\\
\hline
$0.00613945693$	& $0.00626790650$	& $0.00617127881$	& $0.00614739440$	& $2.2151$e$-10$\\
\hline
\end{tabular}
\caption{Record of the $5$ highest eigenvalues of the Karhunen-Loève decomposition of the stationary Ornstein-Uhlenbeck process defined by the SDE $dr_t = -\theta r_t dt + \sigma dW_t, \hspace{2mm} r_0 \sim \mathcal{N}\left(0,\frac{1}{2}\right)$. }
\label{tab:stationary_ornstein_eigenvalues_equi}
}
\end{center}
\end{figure}
With regard to the above numerical results, Nyström method yields a satisfactory accuracy for performing functional quantization of these processes.
\subsection{Eigenfunctions accuracy}
\par \noindent We now compare the closed-form expression of the eigenfunction with the approximation obtained by "Richardson-Romberg extrapolated trapezoidal Nyström method". In table \ref{tab:eigenfunctions_accuracy_tab}, we report the highest absolute difference between the closed-form expression and the approximation on a $300$ points regular mesh of $[0,1]$. The tested cases are the Brownian motion, the Brownian bridge and the stationary Ornstein-Uhlenbeck process defined by the SDE (\ref{eq:sde_ornstein_stationary}) with $\sigma = 1$ and $\theta = 1$.
\begin{figure}[!ht]
\begin{center}
\hspace*{-3mm}
{\small
\begin{tabular}{|c|c|c|c|c|c|}
\hline
Richardson-Romberg	& 						&						&						&						&\\
$50-100-200$  		& {\normalsize $e_1$}	& {\normalsize $e_2$}	& {\normalsize $e_3$}	& {\normalsize $e_4$}	& {\normalsize $e_5$}\\
 absolute error		&						&						&						&						&\\
\hline \hline
Standard			&					&				&				&				&\\
Brownian motion 	& $3.8769$e$-6$		& $3.4909$e$-5$	& $9.6779$e$-5$	& $1.9053$e$-3$	& $3.1558$e$-3$\\
on $[0,1]$			&					&				&				&				&\\
\hline
Standard			&					&				&				&				&\\
Brownian bridge		& $1.5505$e$-5$		& $6.2096$e$-5$	& $1.1398$e$-3$	& $2.4863$e$-3$	& $3.8531$e$-3$\\
 on $[0,1]$ 		&					&				&				&				&\\
\hline
Stationary Ornstein-Uhlenbeck &			&				&				&				&\\
process on $[0,1]$ 	& $3.2257$e$-6$		& $2.1355$e$-5$	& $6.8185$e$-5$	& $1.4614$e$-3$	& $2.5523$e$-3$\\
with $\sigma = 1$ and $\theta = 1$ &	&				&				&				&\\
\hline
\end{tabular}
}
\caption{Record of the biggest absolute error on the Karhunen-Loève eigenfunctions approximation by the Richardson-Romberg extrapolated trapezoidal Nyström method. The number of time steps used for the $3$ steps interpolation are $50$, $100$ and $200$. $300$ equally spaced points on $[0,1]$ were tested. Each column corresponds to one eigenfunction.}
\label{tab:eigenfunctions_accuracy_tab}
\end{center}
\end{figure}

\section{Quantization of the fractional Brownian motion}
\par The normalized fractional Brownian motion $B^H$, is a centered Gaussian process on $[0,T]$, which has the following covariance function:
\begin{equation}\label{eq:fbm_covariance_function}
\Gamma_{B^H}(t,s) = \frac{1}{2}\left( |t|^{2H}+|s|^{2H} - |s-t|^{2H}\right),
\end{equation}
where $H \in(0,1)$ is called the Hurst parameter. If $H = \frac{1}{2}$ then the process is the standard Brownian motion.
\par
\begin{figure}[!ht]
\begin{minipage}[c]{.32\linewidth}
\psfrag{0}{$0$}
\psfrag{0.2}{$0.2$}
\psfrag{0.4}{$0.4$}
\psfrag{0.6}{$0.6$}
\psfrag{0.8}{$0.8$}
\psfrag{1}{$1$}
\psfrag{-2.5}{$-2.5$}
\psfrag{-2}{$-2$}
\psfrag{-1.5}{$-1.5$}
\psfrag{-1}{$-1$}
\psfrag{-0.5}{$-0.5$}
\psfrag{0.5}{$0.5$}
\psfrag{1.5}{$1.5$}
\psfrag{2}{$2$}
\psfrag{2.5}{$2.5$}
\includegraphics[height=4.5cm]{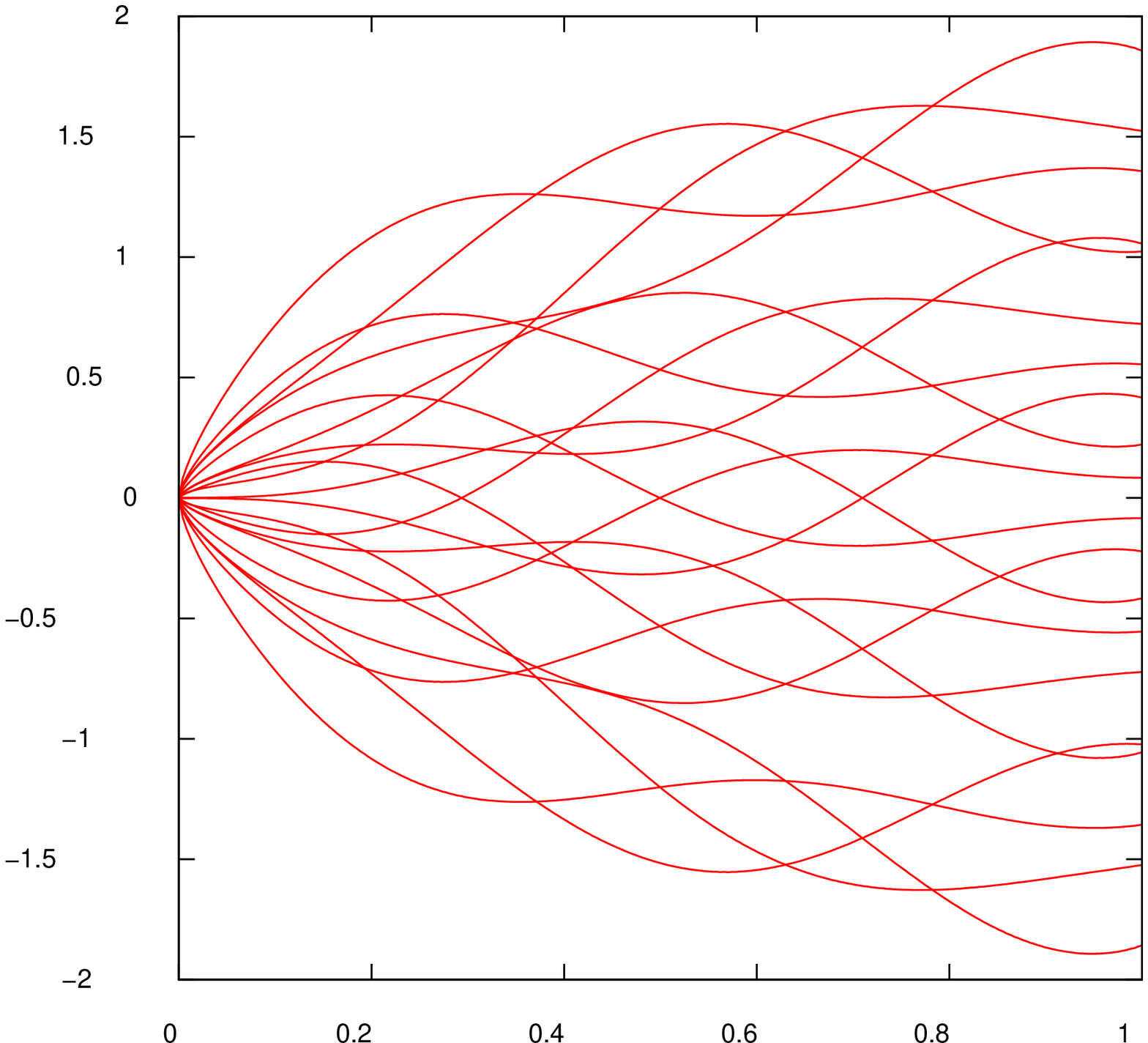}
\end{minipage} \hfill
\begin{minipage}[c]{.32\linewidth}
\psfrag{0}{$0$}
\psfrag{0.2}{$0.2$}
\psfrag{0.4}{$0.4$}
\psfrag{0.6}{$0.6$}
\psfrag{0.8}{$0.8$}
\psfrag{1}{$1$}
\psfrag{-2.5}{$-2.5$}
\psfrag{-2}{$-2$}
\psfrag{-1.5}{$-1.5$}
\psfrag{-1}{$-1$}
\psfrag{-0.5}{$-0.5$}
\psfrag{0.5}{$0.5$}
\psfrag{1.5}{$1.5$}
\psfrag{2}{$2$}
\psfrag{2.5}{$2.5$}
\includegraphics[height=4.5cm]{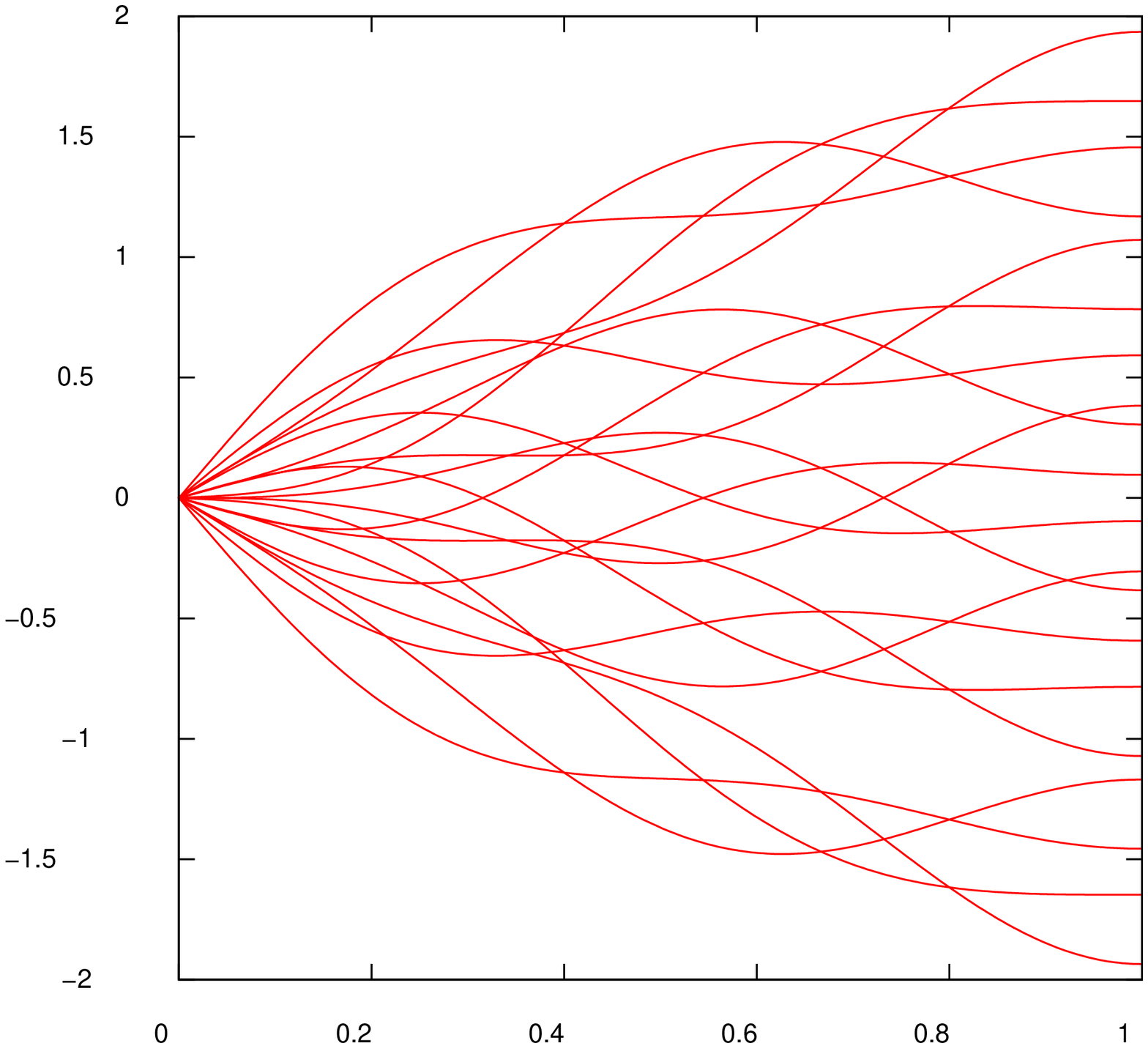}
\end{minipage} \hfill
\begin{minipage}[c]{.32\linewidth}
\psfrag{0}{$0$}
\psfrag{0.2}{$0.2$}
\psfrag{0.4}{$0.4$}
\psfrag{0.6}{$0.6$}
\psfrag{0.8}{$0.8$}
\psfrag{1}{$1$}
\psfrag{-2.5}{$-2.5$}
\psfrag{-2}{$-2$}
\psfrag{-1.5}{$-1.5$}
\psfrag{-1}{$-1$}
\psfrag{-0.5}{$-0.5$}
\psfrag{0.5}{$0.5$}
\psfrag{1.5}{$1.5$}
\psfrag{2}{$2$}
\psfrag{2.5}{$2.5$}
\includegraphics[height=4.5cm]{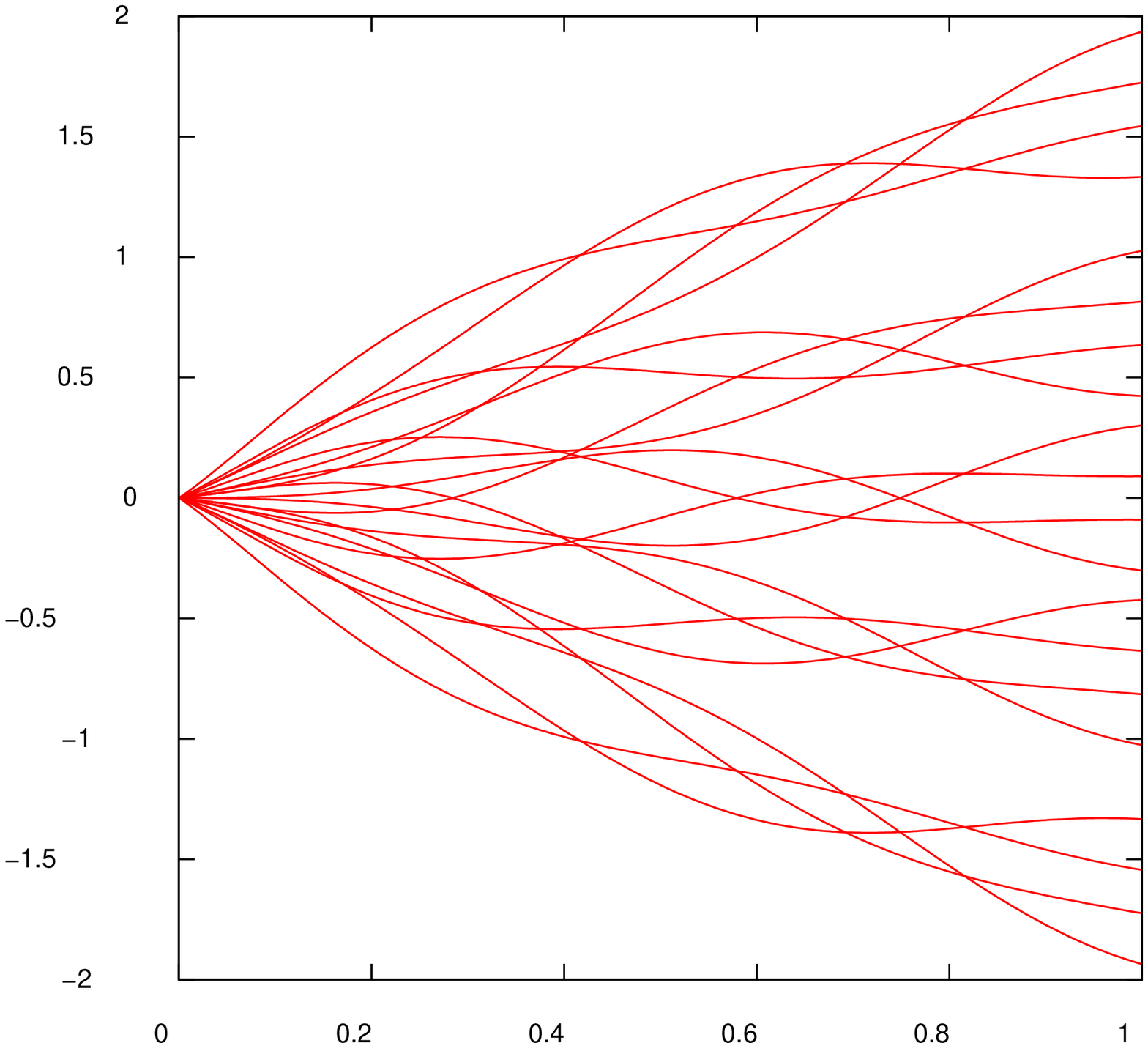}
\end{minipage}
\caption{$(5\times 2 \times 2) -$product quantizer of fractional Brownian motions on $[0,1]$ with Hurst exponent $H=0.3$ (left), $H=0.5$ (middle) and $H=0.7$ (right).}
\label{fig:FBM_product_quantizer}
\end{figure}

 A simple application of the Nyström method presented in section \ref{seq:nystrom_introduction} produces regularly shaped functional quantizers of the fractional Brownian motion. In figure \ref{fig:FBM_product_quantizer}, a $(5\times 2 \times 2) -$product quantizer of the fractional Brownian motion with $3$ different values of the Hurst parameter is plotted.

\par Still, for $H<\frac{1}{2}$, the covariance function of the fractional Brownian motion has singularities that break the convergence of the trapezoidal integration rule in even powers of the step sizes. Indeed, the derivative of $t\to \Gamma_{B^H}(t,s)$ has an infinite limit for $t \to 0^+$ and for ($t \to s^-$ or $t \to s^+$). It breaks also the convergence of the whole associated Nyström method in even powers of the step sizes. In \cite{AtkinsonMonograph,DelvesMohamedIntegral,NumericalRecipes}, methods to handle such boundary and diagonal singularities are proposed. We will deal with this in section \ref{seq:fractional_singularities}
\par However, it is not the case for $H \geq \frac{1}{2}$, so that we can be confident in the results of this method in this case. In table \ref{tab:fractional_07_eigenvalues}, we report the $5$ highest Karhunen-Loève eigenvalues of the fractional Brownian motion on $[0,1]$ with Hurst exponent $H=0.7$. The number of time steps are $128$, $256$ and $512$. Last column yields the corresponding three-steps Richardson-Romberg extrapolation. All the computation has been performed with an octuple precision floating point number implementation to increase the accuracy of the $513 \times 513$-matrix eigensystem computation. (Let us precise that in the case of the Brownian motion on $[0,1]$, when performing the same computation, we get an absolute error smaller than $1$e$-15$ for the five first eigenvalues.)

\begin{figure}[!ht]
\begin{center}
{
\begin{tabular}{|c|c|c|c|}
\hline
Trapezoidal		& Trapezoidal		& Trapezoidal		& Trapezoidal\\
	Nyström		&	Nyström 		&	Nyström 		&	Nyström\\
$128$ points	& $256$ points		& $512$ points 		& $128-256-512$ Richardson-Romberg\\
\hline \hline
$0.374536638$	& $0.374533535$		& $0.374532774$		& $0.374532521757236$\\
\hline
$0.0250351543$	& $0.0250343274$	& $0.0250341354$	& $0.0250340726875501$\\
\hline
$0.00728913038$	& $0.00728860123$	& $0.00728848368$	& $0.0072884458064217$\\
\hline
$0.00322117252$	& $0.00322075790$	& $0.00322066901$	& $0.0032206406932789$\\
\hline
$0.00176153269$	& $0.00176116702$	& $0.00176109039$	& $0.00176106615722872$\\
\hline
\end{tabular}
\caption{Record of the $5$ highest eigenvalues of the fractional Brownian motion on $[0,1]$ with Hurst exponent $H = 0.7$.}
\label{tab:fractional_07_eigenvalues}
}
\end{center}
\end{figure}

\subsection{Kernel singularities when $H < \frac{1}{2}$}\label{seq:fractional_singularities}
\par As pointed out above, the covariance function of the fractional Brownian has a boundary singularity for $t \to 0^+$ and a diagonal singularity. In this section, we will use classical methods to handle this kind of singularities. See \cite{AtkinsonMonograph,DelvesMohamedIntegral,NumericalRecipes} for a review of these method. 
\subsubsection{Handling the boundary singularity}
\par \noindent \textbf{Change of variable} 
\par \noindent The singular behavior of the fractional Brownian motion's covariance function $\Gamma_{B^H}$ defined in equation (\ref{eq:fbm_covariance_function}) can be removed by a change of variable. The change of variable $u = t^{2H}$ and $v = s^{2H}$ in integral (\ref{eq:KL_eigen_equation}) yields:
\begin{equation}\label{eq:fbm_changed}
\int_0^{T^{2H}} \Gamma_{B^H}\left(u^{\frac{1}{2H}},v^{\frac{1}{2H}}\right) f_k\left(v^\frac{1}{2H}\right)\frac{1}{2H} v^{\frac{1}{2H}-1} dv = \lambda_k f_k\left(u^{\frac{1}{2H}}\right). 
\end{equation}
(The second change of variable is done to preserve the symmetry of the Kernel.)
\par \noindent This comes to
\begin{equation}\label{eq:fbm_changed_2}
\int_0^{T^{2H}} \frac{1}{2} \left(|u| + |v|  - |u^{\frac{1}{2H}}-v^\frac{1}{2H}|^{2H} \right) f_k\left(v^\frac{1}{2H}\right) \frac{1}{2H} v^{\frac{1}{2H}-1} dv = \lambda_k f_k \left(u^{\frac{1}{2H}}\right). 
\end{equation}
\par \noindent \textbf{Quadrature rule on a single interval} 
\par We now derive a quadrature rule on $[0,T]$ with respect to the weight function $w(v) = \frac{1}{2H} v^{\frac{1}{2H}-1} = w(v) = \frac{1}{2H} v^{\alpha}$ with $\alpha := \frac{1}{2H}-1$. The aim is to make the quadrature rule exact with affine functions as the trapezoidal quadrature rule is, in the case of an integration with a constant weight. 
$$
\int_{l}^{r} \frac{1}{2H} x^\alpha (a x + b) dx = w_l (a l + b) + w_r (a r + b) \hspace{5mm} \forall (a,b) \in \R^2.
$$
This yields
$$
\frac{1}{2H} \left(\frac{a}{\alpha+2} (r^{\alpha+2} - l^{\alpha+2}) + \frac{b}{\alpha + 1} (r^{\alpha+1} - l^{\alpha+1}) \right) = a (w_l l + w_r r) + b (w_l + w_r) \hspace{5mm} \forall (a,b) \in \R^2.
$$
i.e.
$$
\left(\begin{array}{cccc} l & r \\ 1 & 1 \end{array}\right) \left(\begin{array}{cc} w_l \\ w_r \end{array}\right) = \left(\begin{array}{cc} \frac{1}{2H}\frac{1}{\alpha+2}\left( r^{\alpha + 2} - l^{\alpha + 1}\right) \\ \frac{1}{2H}\frac{1}{\alpha + 1}\left(r^{\alpha +1} - l^{\alpha +1}\right)\end{array}\right).
$$
The solution of the linear system is 
$$
w_l = \frac{1}{2H} \frac{(\alpha+1)l^{\alpha+2} + r^{\alpha+2} - (\alpha+2)l^{\alpha+1}r}{(\alpha +1) (\alpha +2)(r - l)}, \hspace{5mm}
w_r = \frac{1}{2H} \frac{(\alpha+1)r^{\alpha+2} + l^{\alpha+2} - (\alpha+2)r^{\alpha+1}l}{(\alpha +1) (\alpha +2)(r - l)}.
$$
This is
$$
w_l = \frac{ l^{\frac{1}{2H}+1} + 2H r^{\frac{1}{2H}+1} - (2H+1) l^{\frac{1}{2H}}r}{(2H+1)(r - l)}, \hspace{5mm}
w_r = \frac{ r^{\frac{1}{2H}+1} + 2H l^{\frac{1}{2H}+1} - (2H+1) r^{\frac{1}{2H}}l}{(2H+1)(r - l)}.
$$
\par \noindent \textbf{Quadrature rule for equally spaced abscissas} 
\par Let us now consider the equally spaced abscissas points $x_i = i \frac{T}{n}$, $i = 0, 1, \cdots, n$.
We now use these weights $n$ times to integrate on intervals $(x_0^{2H},x_1^{2H}),(x_1^{2H},x_2^{2H}), \cdots, (x_{n-1}^{2H},x_n^{2H})$ to obtain  the extended rule of quadrature. The convergence rate of this method is the same as the trapezoidal rule. 
\subsubsection{Handling the diagonal singularity}
\par We now have to handle the diagonal singularity $\left|u-v \right|^{2H}$ in equation (\ref{eq:KL_eigen_equation}). One classical method if to use the smoothness of the solution by \emph{subtracting of the singularity}.
$$
\int_0^T \Gamma_{B^H} (t,s) f(s) ds = \int_0^{T} \Gamma_{B^H} (t,s) \left( f(s) -f(t) \right) ds + r(t) f(t),
$$
where $r(t) = \int_0^T \Gamma_{B^H} (t,s) ds $. The discretized eigenvalue problem is now transformed to
\begin{equation}\label{eq:diag_eigenvalue}
\begin{array}{lll}
\lambda_k f_k(t_i)	& = \sum\limits_{j=1}^{n} w_j K_{ij} \left(f_k(t_j) - f_k(t_i) \right) + r(t_i) f_k(t_i)\\
					& = \sum\limits_{j=1}^{n} w_j K_{ij} f_k(t_j) + \left(r(t_i) - \sum\limits_{j=0}^{n} w_j K_{ij}\right) f_k(t_i).
\end{array}
\end{equation}
\par \noindent We now define the diagonal matrix $D = \diag(w_i)_{1\leq i \leq n}$ and $D^{1/2} = \diag(\sqrt{w_i})_{1\leq i \leq n}$ as in section \ref{seq:nystrom_introduction}. Moreover, we denote $\Delta = \diag\left(r(t_i) - \sum\limits_{j=0}^{n} w_j K_{ij}\right)_{1\leq i \leq n}$. 
\par \noindent Equation (\ref{eq:diag_eigenvalue}) writes
$$
\lambda_k f_k = K \cdot D f_k + \Delta f_k.
$$
Multiplying by $D^{\frac{1}{2}}$ yields $\lambda h = \left(D^{\frac{1}{2}} \cdot K \cdot D^{\frac{1}{2}} + \Delta \right) h,$ with $h = D^{\frac{1}{2}} f$. As a consequence, we obtain again a symmetric matrix eigenvalue problem. In the case of the fractional Brownian motion, the function $r(t) = \int_0^T \Gamma_{B^H} (t,s) ds$ is derived explicitly:
$$
r(t) = \frac{1}{2} \left( \frac{T^{2H+1}-u^{2H+1}}{2H+1} + u^{2H} T  - \frac{(T-u)^{2H+1}}{2H+1} \right).
$$
\subsubsection{Optimal quantization of the fractional Brownian motion}
\par \noindent We now use this approximation of the Karhunen-Loève basis to perform an optimal quantization of the fractional Brownian motion with a $50$-$100$-$200$ three-step Richardson-Romberg extrapolated Nyström method.
\par
\begin{figure}[!ht]
	\begin{center}
	\psfrag{-2}{$-2$}
	\psfrag{-1.5}{$-1.5$}
	\psfrag{-1}{$-1$}
	\psfrag{-0.5}{$-0.5$}
	\psfrag{0}{$0$}
	\psfrag{0.5}{$0.5$}
	\psfrag{1}{$1$}
	\psfrag{1.5}{$1.5$}
	\psfrag{2}{$2$}
	\psfrag{0.2}{$0.2$}
	\psfrag{0.4}{$0.4$}
	\psfrag{0.6}{$0.6$}
	\psfrag{0.8}{$0.8$}
	\includegraphics[height=6cm]{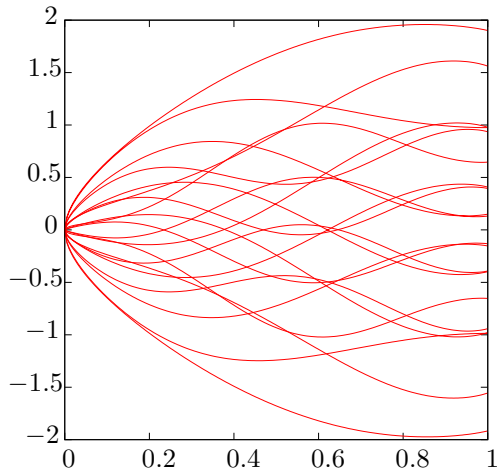}
	\caption{Quadratic $N$-optimal quantizer of the fractional Brownian motion on $[0,1]$ with Hurst's parameter $H=0.25$ and $N=20$.}
	\label{fig:fractional_brownian_motion_optimal_quantization}
	\end{center}
\end{figure}

 In figure \ref{fig:fractional_brownian_motion_optimal_quantization}, we display the quadratic optimal $N-$quantizer of the fractional Brownian motion on $[0,1]$ with Hurst exponent $H = 0.25$ and $N = 20$. In this case, the quantization dimension is $3$. 
\section{Functional stratification of the fractional Brownian motion}
\par In this section, we experiment the functional quantization based stratified sampling algorithm proposed in \cite{CorlayPagesStratification} with the fractional Brownian motion.  
\subsection{Background on stratification}
\par Let $E$ be a separable Hilbert space. The idea of stratification is to localize the Monte-Carlo simulation on the elements of a measurable partition of the state space of a $L^2$ random variable $X:(\Omega,\mathcal{A}) \to (E,\varepsilon)$.
\begin{itemize}
\item Let $(A_i)_{i \in I}$ be a finite $\varepsilon$-measurable partition of a $E$. The sets $A_i$ are called \emph{strata}. Assume that the weights 
$p_i = \PP(X \in A_i)$ are known for $i \in I$ and strictly positive. 
\item Let us define the collection of independent random variables $(X_i)_{i \in I}$ with distribution $\mathcal{L}(X |X \in A_i)$.
\end{itemize}
\par \noindent Let $F:(E,\varepsilon) \to (\R,\mathcal{B}(\R))$ such that $\E[F^2(X)]<+\infty$. 
$$
\begin{array}{lll}
\E[F(X)] 	&= \sum\limits_{i \in I} \E[\1_{\{ X_i\in A_i\}}F(X)] = \sum\limits_{i \in I} p_i \E[F(X) | X \in A_i]\\
			&= \sum\limits_{i \in I} p_i \E[F(X_i)].
\end{array}
$$
\par \noindent The stratification concept comes into play now. Let $M$ be the global budget allocated to the computation of $\E[F(X)]$ and $M_i=q_i M$ the budget allocated to compute $\E[F(X_i)]$ in each stratus. We assume that $\sum\limits_{i \in I} q_i = 1$. This leads to define the (unbiased) estimator of $\E[F(X)]$:
\begin{equation}\label{eq:unbiased_stratif_estimator}
\overline{F(X)}^I_M := \sum\limits_{i \in I} p_i \frac{1}{M_i} \sum\limits_{k=1}^{M_i} F(X_i^k),
\end{equation}
where $(X_i^k)_{1\leq k \leq M_i}$ is a $\mathcal{L}(X|X \in A_i)$-distributed random sample. 
\begin{prop}
\par With the same notations:
\begin{equation}\label{eq:stratif_estimator_variance}
\var\big(\overline{F(X)}^I_M\big)=\frac{1}{M} \sum\limits_{i\in I} \frac{p_i^2}{q_i} \sigma_{F,i}^2,
\end{equation}
where $\sigma^2_{F,i} = \var(F(X)|X\in A_i) = \var(F(X_i)) \ \forall i \in I$.
\end{prop}
\par \noindent The proof can be found in \cite{CorlayPagesStratification}. Optimizing the simulation allocation to each stratus amounts to solving the following minimization problem: 
\begin{equation}\label{eq:stratification_minimization}
\min\limits_{(q_i) \in \mathcal{P}_I} \sum\limits_{i \in I} \frac{p_i^2}{q_i} \sigma^2_{F,i} \hspace{2mm} \textrm{ where } \mathcal{P}_I = \left\{ (q_i)_{i\in I} \in \R_+^I \Big| \sum\limits_{i\in I} q_i = 1 \right\}. 
\end{equation}
\par \noindent In \cite{CorlayPagesStratification}, Corlay and Pagès pointed out theoretical aspects of quantization that lead to a strong link between the problem of optimal $L^2$-quantization of a random variable and the variance reduction that can be achieved by stratification. Three types of allocation rules for the budgets $(q_i)_{i \in I}$ are proposed:
\begin{itemize}
\item The "sub-optimal rule" is to set
\begin{equation}\label{eq:sub_optimal_choice}
q_i=p_i, \hspace{5mm} i \in I.
\end{equation}
The two motivations for this choice are the facts that the weights $p_i$ are known and because it always reduces the variance. 
\item The "optimal rule" is the solution of the constrained minimization problem (\ref{eq:stratification_minimization}). The Schwartz inequality yields
$$
\sum\limits_{i \in I} p_i \sigma_{F,i} = \sum\limits_{i \in I} \frac{p_i \sigma_{F,i}}{\sqrt{q_i}} \sqrt{q_i} \leq \left(\sum\limits_{i \in I} \frac{p_i^2 \sigma_{F,i}^2}{q_i} \right)^{1/2}{\underbrace{\left( \sum\limits_{i\in I} q_i \right)}_{=1}}^{1/2}.
$$
\par As a consequence, the solution of the minimization problem corresponds to the equality case into the Schwartz inequality. Hence the solution of the minimization problem is given by 
\begin{equation}\label{eq:optimal_budget}
q_i^* = \frac{p_i \sigma_{F,i}}{\sum\limits_{j\in I} p_j \sigma_{F,j}}, i\in I
\end{equation}
and the corresponding minimal variance is given by $\left(\sum\limits_{i \in I} p_i \sigma_{F,i}\right)^2.$
\par The counterpart of this method is that we do not know explicitly the solution $(q_i^*)_{i \in I}$. In \cite{JourdainStratification1}, \'{E}toré and Jourdain proposed an algorithm for adaptively modifying the proportion of further drawings in each stratum, that converges to the optimal allocation. This can be used in a general framework. Another practical solution would be to implement a simple prior rough estimation of the optimal allocation.
\item The "Lipschitz optimal" rule. When the partition $(A_i)_{i \in I}$ is a Voronoi partition associated with an optimal quantizer of $X$, Corlay and Pagès considered the setting
\begin{equation}\label{eq:universal_stratification}
q_i=\sigma_i, \hspace{5mm} i \in I,
\end{equation}
where $\sigma_i$ is the local inertia of the random variable $X$, $\sigma_i^2 = \E\Big[|X-\E[X|X\in A_i]|^2 \Big| X \in A_i \Big].$ It is proved that this setting has a uniform efficiency among the class of Lipschitz continuous functionals. Moreover, local inertia $(\sigma_i)_{i \in I}$ are known. This solution overcomes the "sub-optimal choice" in every test done in \cite{CorlayPagesStratification}.
\end{itemize}
\subsection{On the functional stratification of Gaussian processes}
\par Here, we assume that $X$ is an $\R$-valued Gaussian process on $[0,T]$. We are interested in the value of $\E[F(X_{t_0}, X_{t_1}, \cdots, X_{t_n})]$ where $0 = t_0 \leq t_1 \leq \cdots \leq t_n = T$ are $n+1$ dates of interest for the underlying process. Let us assume that $\chi \in \mathcal{O}_{pq}(X,N)$ is a K-L product quantizer of $X$. The codebook associated with this product quantizer is the set of the paths of the form 
$$
\chi_{\ui} = \sum\limits_{n \geq 1} \sqrt{\lambda_n^X} x_{i_n}^{(N_n)} e_n^X, \hspace{5mm} \ui = \{i_1, \cdots, i_n, \cdots \},
$$
where $(e_n^X,\lambda_n^X)$ is the Karhunen-Loève decomposition of the process $X$ on $[0,T]$ and $x_{i_n}^{N_n}$ is the $i_n$th element of an optimal quantizer of size $N_n$ of the standard one-dimensional Gaussian distribution. 
\par We now need to be able to simulate the conditional distribution 
$$
\mathcal{L}(X | X \in A_{\ui})
$$
where $A_{\ui}$ is the slab associated with $\chi_{\ui}$ in the codebook. 
\par To simulate the conditional distribution $\mathcal{L}(X | X \in A_{\ui})$, we will:
\begin{itemize}
	\item First, simulate the first K-L coordinates of $X$. The explicit simulation algorithm is available in \cite{CorlayPagesStratification}
	\item Then simulate the conditional distribution of the marginals of the Gaussian process, its first coordinates being settled. 
\end{itemize}
\par \noindent In this setting, the aim is to simulate the conditional distribution 
\begin{equation}\label{eq:conditional_distribution}
\mathcal{L}\Big(X_{t_0},\cdots,X_{t_n} \Big| \int_0^T X_s e_1^X ds, \int_0^T X_s e_2^X(s) ds, \cdots , \int_0^T X_s e_d^X(s) ds \Big)
\end{equation}
where $(X_t)_{t \in} [0,T]$ is a $L^2$ $\R$-valued Gaussian process, and $(e_k^X,\lambda_k^X)_{k\in \N^*}$ is the Karhunen-Loève system associated with the process $X$.
\par \noindent \textbf{Conditional simulation:} In \cite{CorlayPagesStratification}, two solutions are proposed for the simulation of the conditional distribution (\ref{eq:conditional_distribution}). 
\begin{itemize}
\item The first one is the naive Cholesky method for Gaussian vector simulation, which has a quadratic complexity in the number of time steps. This first simulation scheme was not competitive for linearly simulable processes as the Brownian motion. In the following, we will mention this method as the \emph{brute force method}.
\item The other solution, detailed in \cite{CorlayPagesStratification} requires a prior simulation of the unconditional distribution of $(X_{t_0},\cdots,X_{t_n})$ and has then a linear additional cost. This algorithm will be mentioned in the following as the \emph{linear conditioning algorithm}. For Gaussian processes which have a linear simulation scheme in the unconditional case (as the Ornstein-Uhlenbeck process, the Brownian bridge and the Brownian motion), this method is of high interest.
\end{itemize}
\subsection{The case of the fractional Brownian motion}
\par Possible methods for simulating the fractional Brownian motion on a schedule $t_0 < t_1 < \cdots < t_n$ are 
\begin{itemize}
\item the naive Cholesky method, that has quadratic complexity,
\item and the circulant matrix method which has a $O(n \ln (n))$ complexity \cite{DietrichSimulationFractional,WoodChanSimulationStationary}. The circulant matrix method is also available for the multifractional Brownian motion \cite{WoodChanSimulationMultifractional}.
\end{itemize}
No exact simulation scheme with a linear complexity exists for the fractional Brownian motion. Still, approximate method with linear complexity exists. If we choose the Cholesky method, there is no interest to use the linear conditioning algorithm proposed in \cite{CorlayPagesStratification}. The brute force Cholesky method is adapted to this situation.
\par In every other case, if the unconditional simulation method has smaller complexity, we have interest to use the linear conditioning algorithm which has a linear additional cost to the unconditional simulation. 
\par
\begin{figure}[!ht]
	\begin{center}
	\psfrag{-4}{$-4$}
	\psfrag{-3}{$-3$}
	\psfrag{-2}{$-2$}
	\psfrag{-1}{$-1$}
	\psfrag{0}{$0$}
	\psfrag{1}{$1$}
	\psfrag{2}{$2$}
	\psfrag{3}{$3$}
	\psfrag{4}{$4$}
	\psfrag{0.5}{$0.5$}
	\psfrag{1.5}{$1.5$}
	\psfrag{2.5}{$2.5$}
	\includegraphics[height=6cm]{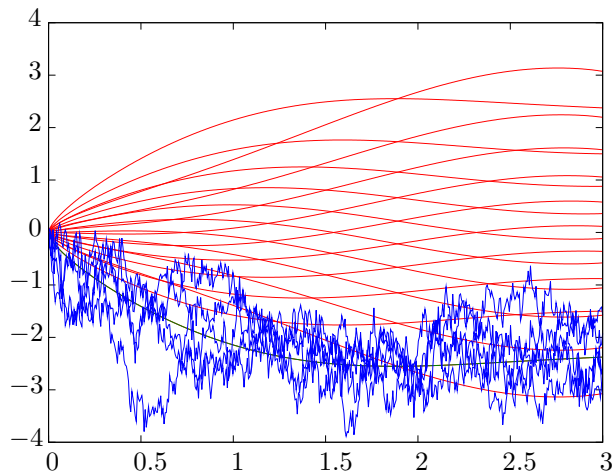}
	\caption{Plot of a few paths of the conditional distribution of the fractional Brownian motion with Hurst's parameter $H=0.3$ on $[0,3]$, knowing that its path belong to the $L^2$ Voronoi cell of the highlighted curve in the quantizer. }
	\label{fig:fractional_brownian_motion_conditional_simulation}
	\end{center}
\end{figure}

 In figure \ref{fig:fractional_brownian_motion_conditional_simulation}, we plot a few paths of the conditional distribution of the fractional Brownian motion with Hurst's parameter $H=0.3$ knowing that they belong to a given $L^2$ Voronoi cell. 
\subsection{Gaussian process reconstruction}
\par The first numerical test of the functional stratification of the fractional Brownian motion is a method to validate both the eigenfunction computation by the Nyström method and the functional stratification algorithm. 
\par Indeed, one can rebuild the considered Gaussian process from its stratification. This yields the following simulation algorithm:
\begin{itemize}
\item First, simulate the discrete weighted distribution of the strata index $(i,p_i)_{i \in I}$ to select the strata.
\item Then simulate the conditional distribution $\mathcal{L}\left(X_{t_0}, \cdots, X_{t_n} \Big|X\in A_i \right)$  of the Gaussian process in the strata by the method described above.
\end{itemize}
The result should be distributed according to the distribution of the underlying Gaussian process. In table \ref{tab:fractional_reconstruction_07}, we report the covariance structure $\E[X_{t_i}X_{t_j}]_{1 \leq i,j \leq n}$ estimated by a Monte-Carlo simulation when $X$ is a fractional Brownian motion with Hurst's parameter $H=0.7$. The tested schedule is $(i \frac{T}{n})_{0 \leq i \leq n}$ with $T = 1$ and $n = 5$. The product decomposition of the quantization is $10 \times 5 \times 2$.
\begin{figure}[!ht]
{\small
\hspace*{-7mm}
\begin{minipage}[c]{.46\linewidth}
\begin{center}
\begin{tabular}{|c|c|c|c|c|}
\hline
$0.105061$ & $0.138629$ & $0.15846$ & $0.173817$ & $0.186687$\\
\hline
$0.138629$ & $0.277258$ & $0.330656$ & $0.365844$ & $0.394071$\\
\hline
$0.15846$ & $0.330656$ & $0.489116$ & $0.557871$ & $0.605929$\\
\hline
$0.173817$ & $0.365844$ & $0.557871$ & $0.73168$ & $0.813313$\\
\hline
$0.186687$ & $0.394071$ & $0.605929$ & $0.813313$ & $1$\\
\hline
\end{tabular}
\end{center}
\end{minipage} \hfill
\begin{minipage}[c]{.46\linewidth}
\begin{center}
\begin{tabular}{|c|c|c|c|c|}
\hline
$0.105141$ & $0.138748$ & $0.158596$ & $0.173959$ & $0.186824$\\
\hline
$0.138748$ & $0.277417$ & $0.330885$ & $0.366075$  & $0.394372$\\
\hline
$0.158596$ & $0.330885$ & $0.489454$ & $0.558177$ & $0.606266$\\
\hline
$0.173959$ & $0.366075$ & $0.558177$ & $0.731923$ & $0.813579$\\
\hline
$0.186824$ & $0.394372$ & $0.606266$ & $0.813579$  & $1.0003$\\
\hline
\end{tabular}
\end{center}
\end{minipage}
}
\caption{Theoretical (left) and estimated (right) covariance $\E[X_{t_i} X_{t_j}]$ of the rebuilt fractional Brownian motion with $H = 0.7$. The number of generated paths for this Monte-Carlo simulation was $1\cdot 10^7$.}
\label{tab:fractional_reconstruction_07}
\end{figure}
\par \noindent In every tested case, when generating table \ref{tab:fractional_reconstruction_07}, the theoretical value lies in the $95\%$ confidence interval. These confidence intervals were not displayed for briefness. We obtain the same order of accuracy with other values of $H \in (0,1)$.
\subsection{Application to option pricing}
\par A stochastic integral with respect to the fractional Brownian motion has been introduced in \cite{HeliotHoekFractional} by Helliot and van der Hoek, and in \cite{OksendalWhiteNoise} by Biagini, Øksendal, Sulem and Wallner. They proposed a generalization of the Black-Scholes model. As in the classical Black-Scholes market, two assets are available: 
\begin{itemize}
\item A risk-free asset whose price is given by
\begin{equation}\label{eq:non_risky_BS}
dS^0_t = r S^0_t dt
\end{equation}
\item and a risky asset whose price is given by
\begin{equation}\label{eq:fractional_BS_SDE}
d S_t = \mu S_t dt + \sigma S_t d B_t^H,
\end{equation}
where $r$, $\mu$ and $\sigma$ are constants and $B^H$ is fractional Brownian motion with Hurst parameter $H$.
\end{itemize}
\par\noindent It has been shown that this market presents no arbitrage opportunity and is complete. Moreover, the solution of the stochastic differential equation (\ref{eq:fractional_BS_SDE}) is given by
\begin{equation}\label{eq:fractional_BS_SDE_solution}
S_t = S_0 \exp\left(\sigma B^H_t + \mu t - \frac{1}{2}\sigma^2 t^{2H} \right).
\end{equation}
The following theorem, prooved in \cite{HeliotHoekFractional} deals with the price of a European call option.
\begin{theo}[Fractional Black-Scholes Formula]
The price at every time $t \in [0,T]$ of a European call option with strike price $K$ and maturity $T$ is given by
\begin{equation}
C(t,S_t) = S_t \mathcal{N}(d_1)- K e^{-r(T-y)}\mathcal{N}(d_2)
\end{equation}
where
\begin{equation}\label{eq:fractional_d1}
d_1 = \frac{\ln\left( \frac{S_t}{K}\right) + r(T-t) + \frac{\sigma^2}{2} (T^{2H}-t^{2H})}{\sigma \sqrt{T^{2H} - t^{2H}}}
\end{equation}
\begin{equation}\label{eq:fractional_d2}
d_2 = \frac{\ln\left( \frac{S_t}{K}\right) + r(T-t) - \frac{\sigma^2}{2} (T^{2H}-t^{2H})}{\sigma \sqrt{T^{2H} - t^{2H}}}
\end{equation}
\end{theo}
\par This closed-form expression is used to benchmark our simulation scheme of the fractional Brownian motion. 
\subsubsection{Benchmark with a Barrier option in the fractional Black and Scholes model}
\par Here, we benchmark the numerical method for a path dependent option in the case of a Barrier option in the fractional Black and Scholes model. For the sake of simplicity, we consider a log-normal Black and Scholes diffusion with no drift (no interest rate and no dividend). The chosen Hurst exponent is $H=0.3$. The numerical results are reported in table \ref{fig:var_reduc_100_points}.
\par The results are displayed for different values of the initial spot $S$, the strike $K$, the barrier $B$, the volatility $\sigma$, the maturity $T$ and the number of equally spaced fixing dates $n$. 
\par In this table, the first column corresponds to a simple Monte-Carlo estimator. The last three columns correspond to a stratified sampling estimator with different simulation allocation for each strata. 
\par The "sub-optimal weights" column stands for the allocation budget of equation (\ref{eq:sub_optimal_choice}). The "Lip.-optimal weights" column stand for the "universal stratification" budget allocation of equation (\ref{eq:universal_stratification}). Both these two case have explicit allocation rules. Last column, "Optimal weights" corresponds to an estimation of the optimal budget allocation given in expression (\ref{eq:optimal_budget}). 
\begin{figure}[!ht]
	\begin{center}
	{
	%\hspace*{-13mm}
	\scriptsize
	\begin{tabular}{|c|c|c|c|c|c|}
	\hline
							& Simple 				&	Strat. Estimator	& Strat. Estimator		& Strat. Estimator		\tabularnewline
	Parameters				& Estimator				&	sub-optimal weights	& Lip.-optimal weights	& Optimal weights		\tabularnewline
							& 						&						& 						&						\tabularnewline
	\hline
	\hline
	$S = 100$, $K = 100$	& $12.5947$				& $12.5674$				& $12.5566$				& $12.5890$				\tabularnewline	
	$B=125$, $\sigma=0.3$,	& $[12.4429,12.7466]$	& $[12.4732,12.6615]$	& $[12.4654,12.6477]$	& $[12.5201,12.6579]$	\tabularnewline
	$T=1.5$, $n=11$ 		& $\var = 600.5711$		& $\var = 230.8692$		& $\var = 216.3442$		& $\var = 123.5426$		\tabularnewline
	\hline
	$S = 100$, $K = 100$	& $1.3412$				& $1.3826$				& $1.3613$				& $1.3769$ 				\tabularnewline	
	$B=200$, $\sigma=0.3$,	& $[1.2677,1.4146]$		& $[1.3140,1.4511]$		& $[1.3002,1.4224]$		& $[1.3530,1.4009]$		\tabularnewline
	$T=1$, $n=11$ 			& $\var = 140.5978$		& $\var = 122.2808$		& $\var = 97.1538$		& $\var = 14.9352$		\tabularnewline
	\hline
	\end{tabular}
	}
	\caption{Numerical results for the Up In Call option, with $100 = \times 5 \times 2 \time 2$ stratas. }
	\label{fig:var_reduc_100_points}
	\end{center}
\end{figure}
\par \noindent We notice that the quantization based stratified sampling method reduces noticeably the variance of the Monte-Carlo estimator. The universal stratification allocation rule (\ref{eq:universal_stratification}) proposed in \cite{CorlayPagesStratification} overcomes the sub-optimal weight allocation. Moreover, the "optimal allocation" estimation yields a better variance reduction factor. 
\bibliography{biblio}
\end{document}